\documentclass[10pt]{amsart}

\usepackage[leqno]{amsmath}
\usepackage{amssymb}
\usepackage{array}
\usepackage{rotating}
\usepackage{amsthm}

\newtheorem*{lemX}{Lemma}\newtheorem*{thmX}{Theorem}

\newcommand{\set}[1]{\left\{#1\right\}}
\newcommand{\abs}[1]{\left\vert#1\right\vert}
\newcommand{\Res}{\mathrm{Res}}
\begin{document}
\title{Implicit higher derivatives, and a formula of Comtet and Fiolet}\author{Tom Wilde}
\begin{abstract}
Let $F(x,y)$ be a function of two variables, and suppose $y=f(x)$ satisfies $F(x,y)=0$ for $x$ in some range. L. Comtet and M. Fiolet (\cite{Comtet1974}) stated a combinatorial formula for the derivatives $d^ny/dx^n$ in terms of the partial derivatives of $F.$ Their formula, however, contains errors. Given the basic nature of this problem, it seems worthwhile to give the corrected formula, which is in fact slightly simpler than as stated in \cite{Comtet1974}. We give two derivations of the formula, by Lagrange inversion and by induction. We also correct Comtet and Fiolet's expression in \cite{Comtet1974} for the number of terms in the formula.
\end{abstract}
\maketitle

\section{Introduction}\label{S1}
Let $F(x,y)$ be a smooth function of two variables, and suppose $y=f(x)$ is a smooth function satisfying $F(x,y)=0$ for $x$ in some range. It follows that $dF=F_xdx+F_y\frac{dy}{dx}=0$ where $F_x,F_y$ are the partial derivatives of $F$ with respect to $x$ and $y.$ Hence
\begin{equation}\label{dydx}\frac{dy}{dx}=-\frac{F_x}{F_y},\end{equation} as is well known. We can calculate the higher derivatives $d^ny/dx^n,n\geq 2$ by taking the total derivative of equation
(\ref{dydx}), giving 
\begin{equation}\label{TD}
\frac{d^ny}{dx^n}=[\frac{\partial}{\partial x}+\frac{dy}{dx}\frac{\partial}{\partial y}](\frac{d^{n-1}y}{dx^{n-1}})=[\frac{\partial}{\partial x}-\frac{F_x}{F_y}\frac{\partial}{\partial y}]^{n-1}(-\frac{F_x}{F_y}).
\end{equation}
For any given $n,$ this can be expanded to give an expression, analogous to equation (\ref{dydx}), for $d^ny/dx^n$ in terms of partial derivatives of $F,$ but it is not clear how to give the resulting expressions as a uniform formula, valid for all $n.$ L. Comtet and M. Fiolet stated such a formula in $1974$ (\cite[Th\'eor\`eme 1]{Comtet1974}, or see equations (\ref{CF1})-(\ref{CF2}) below), but their formula is not correct. Comtet and Fiolet indicated a proof that is undoubtedly correct in principle, using Lagrange inversion, but did not give their detailed calculations, so we do not know how the errors arose. 

The purpose of this note to give the correct result, which we do at equation (\ref{main}) below. We also give two derivations, one using a form of Lagrange inversion (see equation (\ref{GLI1}) below), and another, elementary proof using induction.

This note is organised as follows. In Section \ref{S2}, we give some notation and in Section \ref{S3}, we give our corrected version of Comtet and Fiolet's formula. A derivation using a Lagrange inversion method is in Section \ref{S4}, and another by induction is given in Section \ref{S5}. In Section \ref{S6}, we discuss Comtet and Fiolet's original formula. In Section \ref{S7}, we observe that Comtet and Fiolet's expression \cite[Th\'eor\`eme 2]{Comtet1974} for the number of terms in their formula, is also incorrect, and we give the correct expression.

\section{Notation: Two-dimensional partitions}\label{S2}
To state our version of Comtet and Fiolet's formula, we need to introduce the concepts of one- and two-dimensional partitions, and set up some notation. A partition of the positive integer $n,$ is a sequence of positive integers $p_1\geq p_2\geq ...\geq p_r>0$ such that $\sum_ip_i=n.$ We write $p\vdash n$ to indicate that $p=(p_i)$ is a partition of $n.$ The $p_i$ are called the parts of the partition, and the number $r$ of parts is denoted $\abs p.$

Analogously, if $m$ and $n$ are non-negative integers, not both zero, then a partition of $(m,n)$ is a lexicographically ordered sequence
\begin{equation}
(p_1,q_1)\geq (p_2,q_2)...\geq (p_r,q_r)>0,
\end{equation} 
where for each $i,$ $p_i$ and $q_i$ are non-negative integers and not both zero, and $\sum_i p_i=m$ and $\sum_i q_i=n.$ The lexicographic ordering is given by the rule $(p,q)>(p^\prime,q^\prime)$ if $p>p^\prime$ or $p=p^\prime$ and $q>q^\prime.$ We write $p\vdash (m,n)$ to denote that $p$ is a partition of $(m,n)$ in the above sense. Confusion will not arise, since it is clear from what comes after the $\vdash,$ whether a one- or two- dimensional partition is indicated. Also as in the one-dimensional case, if $p\vdash (m,n)$ then we write $\abs p$ to denote the number $r$ of terms in the sum. 

For an ordinary partition $p=(p_i)\vdash n,$ let $e_{p,k}$ denote the number of times the number $k$ appears among the $p_i.$ Thus $e_{p,k}=\abs{\set{i\mid p_i=k}}.$ Naturally, for $p\vdash (m,n),$ we set $e_{p,k,l}=\abs{\set{i\mid p_{i1}=k\text{ and }p_{i2}=l}},$ where $p_i=(p_{i1},p_{i2}).$ That is, $e_{p,k,l}$ is the number of times $(k,l)$ appears in the partition. 

Next, we will encounter certain numbers $\alpha_p$ where for an ordinary partition $p\vdash n,$ $\alpha_p$ is defined by
\begin{equation}\label{alpha1}
\alpha_p=\frac{n!}{\prod_i p_i!\prod_k e_{p,k}!},
\end{equation}
and in the two dimensional case $p\vdash (m,n),$ by an obvious analogy,
\begin{equation}\label{alpha2}
\alpha_p=\frac{n!m!}{\prod_i p_{i1}!p_{i2}!\prod_{k,l} e_{p,k,l}!},
\end{equation}
where again the $i^\mathrm{th}$ term of $p$ is being written as $(p_{i1},p_{i2}).$ The coefficients $\alpha_p$ in equation (\ref{alpha1}) appear in Fa\`a di Bruno's well known formula, which we will use in our derivation (see equation \ref{FdB1}), while those in equation (\ref{alpha2}) appear in our main result at equation (\ref{main}). It seems noteworthy that the coefficients of these two formulae are analogous in this way. Finally, for a function $g$ of one or two variables, we use the convenient notation, 
\begin{equation}
g_{x,p}=\prod_i\frac{\partial^{p_i}g}{\partial x^{p_i}}\text{ and }g_{x,y,p}=\prod_i\frac{\partial^{p_{i1}+p_{i2}}g}{\partial x^{p_{i1}}\partial y^{p_{i2}}}.\end{equation}
Again, note that it is clear from the context, whether $p$ is an ordinary or a two dimensional partition. We continue to use $g_x$ and $g_y$ to denote the ordinary partial derivatives $\partial g/\partial x$ and $\partial g/\partial y$ respectively.

\section{The formula}\label{S3}
We now state our corrected version of Comtet's and Fiolet's formula. Let $F(x,y)$ be a smooth function of two variables, and suppose $y=f(x)$ is a smooth function satisfying $F(x,y)=0$ for $x$ in some range. Then at a point where $F_y\neq 0,$ we have
\begin{equation}\boxed{\label{main}
\frac{d^ny}{dx^n}=\sum_{p:\tiny{\begin{cases}p\vdash (n,\abs{p}-1)\\(0,1)\notin p\end{cases}}}(-1)^{\abs p}\alpha_p\frac{F_{x,y,p}}{F_y^{\abs p}}.
}\end{equation}
The notation is as set out in Section 2; in particular, the coefficients $\alpha_p$ are given by equation (\ref{alpha2}). See Section \ref{S6} for a discussion of where Comtet and Fiolet's formula differs from equation (\ref{main}). The sum in equation (\ref{main}) is over two-dimensional partitions $p\vdash (n,r-1)$ where $\abs{p}=r,$ and $(0,1)$ cannot be a part of the partition. 

For example, taking $n=2,$ by inspection there are three partitions $p\vdash (2,\abs{p}-1)$ with $(0,1)\notin p,$ namely 
\begin{equation}
(2,0),(1,1)+(1,0)\text{ and }(1,0)+(1,0)+(0,2),
\end{equation} 
and these have $\alpha_p=1,2$ and $1$ respectively, so we obtain
\begin{equation}\label{EG}
\frac{d^2y}{dx^2}=-\frac{F_{xx}}{F_y}+2\frac{F_xF_{xy}}{F_y^2}-\frac{F_x^2F_{yy}}{F_y^3},
\end{equation}
as may easily be verified using equation (\ref{TD}).

\section{Derivation of equation (\ref{main})}\label{S4}
Comtet and Fiolet's indicated proof of their equation refers to an earlier paper \cite{Comtet1968}, in which a further indication is given of a proof using Lagrange inversion. However, it seems more direct to use the following observation. In equation (\ref{dydx}), provided $F_y\neq 0,$ $F_x/F_y$ is the residue of $F_x/F$ at its pole in the $y$-plane at $y=f(x).$ Therefore,
\begin{equation}
\frac{dy}{dx}=-\Res_{y=f(x)}(\frac{F_x}{F})=-\frac{1}{2\pi i}\int_{C} dy \frac{F_x}{F}
\end{equation}
where $C$ is a small loop around the pole
at $y=f(x).$ The location of the pole varies with $x,$ but for sufficiently small variation of $x,$ it remains inside $C,$ so we may
differentiate under the integral sign without moving $C.$ Hence
\begin{equation}
\frac{d^ny}{dx^n}=-\frac{1}{2\pi i}\int_C dy \frac{\partial^{n-1}}{\partial x^{n-1}}(\frac{F_x}{F})
=-\Res(\frac{\partial^{n-1}}{\partial x^{n-1}}(\frac{F_x}{F})),
\end{equation} 
where the residue is at $y=f(x)$ in the $y$-plane. We write this concisely as
\begin{equation}
\label{GLI1}\frac{d^ny}{dx^n}=-\Res_y({\frac{\partial^n\log F}{\partial x^n}}),
\end{equation}
where the subscript $y$ refers to the residue being in the $y$ plane.

Equation (\ref{GLI1}) is a relative of Lagrange inversion. The classic Lagrange inversion formula can be recovered by letting $F(x,y)=x-g(y),$ so that $f$ and $g$ are inverse, and evaluating the residue. However, equation (\ref{GLI1}) is particularly suited to our present problem; all we have to do is evaluate the residue on the right hand side. To do this, we first introduce Fa\`{a} di Bruno's formula. If $z$ is a sufficiently differentiable function of $y,$ and $y$ is a sufficiently differentiable function of $x,$ then Fa\`a di Bruno's formula states that
\begin{equation}\label{FdB1}
\frac{d^nz}{dx^n}=\sum_{p\vdash n}\alpha_p\frac{d^{\abs p}z}{dy^{\abs p}}y_{x,p},
\end{equation} 
where the coefficients $\alpha_p$ are given by equation (\ref{alpha1}). The formula can be proved easily by induction. See \cite[Section 24.1.2(C)]{AbSteg} or \cite{Johnson} for details and background.

We can now begin our derivation of equation (\ref{main}). Applying Fa\`{a} di Bruno's formula (\ref{FdB1}) to $\log F(x,y)$ as a function of $x,$ we have
\begin{equation}
\frac{\partial^n}{\partial x^n}\log F=\sum_{p\vdash n}\alpha_p(-1)^{{\abs{p}}-1}({\abs{p}}-1)!\frac{F_{x,p}}{F^{\abs{p}}}.
\end{equation} 
Hence equation (\ref{GLI1}) becomes
\begin{equation}\label{GLI2}
\frac{d^ny}{dx^n}=\sum_{p\vdash n}\alpha_p(-1)^{\abs{p}}({\abs{p}}-1)!\Res_y(\frac{F_{x,p}}{F^{\abs{p}}}).
\end{equation}
Now to simplify the notation, we temporarily work at the origin $x=y=0,$ assuming therefore that $F(0,0)=0$ and $F_y(0)\neq 0.$ Then the pole of
$F_{x,p}/F^{\abs{p}}$ at $y=0$ has order $\abs{p},$ so
\begin{equation}\label{LR0}
\Res_{y=0}(\frac{F_{x,p}}{F^{\abs{p}}})=\frac{1}{(\abs{p}-1)!}\frac{\partial^{\abs{p}-1}}{\partial
y^{\abs{p}-1}}(F_{x,p}(F/y)^{-\abs{p}})\vert_{y=0}.
\end{equation}
By Leibniz' rule,
\begin{equation}\label{LR1}
\frac{\partial^{\abs{p}-1}}{\partial y^{\abs{p}-1}}(F_{x,p}(F/y)^{-\abs{p}})=
\sum_{s=0}^{\abs{p}-1}\frac{(\abs{p}-1)!}{s!(\abs{p}-1-s)!}\frac{\partial^{\abs{p}-1-s}F_{x,p}}{\partial
y^{\abs{p}-1-s}}\frac{\partial^s}{\partial y^s}((F/y)^{-\abs{p}}),
\end{equation}
and by Fa\`{a} di Bruno's formula (\ref{FdB1}) applied to $(F/y)^{-\abs{p}}$ as a function of $F/y,$ we have
\begin{equation}\label{LR2}
\frac{\partial^s}{\partial y^s}((F/y)^{-\abs{p}})=\sum_{q\vdash s}\alpha_q(-1)^{\abs{q}}\frac{(\abs{p}+\abs{q}-1)!}
{(\abs{p}-1)!}\frac{(F/y)_{y,q}}{(F/y)^{\abs{p}+\abs{q}}}.
\end{equation}
To take the limit as $y$ tends to $0,$ in equation (\ref{LR2}), we replace $\partial^r(F/y)/\partial y^r $ with $\frac{1}{r+1}\partial^{r+1}F/\partial y^{r+1}$ for $r\geq 0.$ Hence $F/y$ becomes $F_y$ and 
\begin{equation}(F/y)_{y,q}=\prod_{i\geq 1}(\frac{\partial^{q_i}(F/y)}{\partial y^{q_i}})\text{ becomes } \frac{1}{\prod_{i\geq 1}(q_i+1)}(F_y)_{y,q}.
\end{equation}
Therefore at $y=0,$ equation (\ref{LR2}) becomes
\begin{equation}
\frac{\partial^s}{\partial y^s}((F/y)^{-\abs{p}})=\sum_{q\vdash s}\alpha_q(-1)^{\abs{q}}\frac{(\abs{p}+\abs{q}-1)!}
{(\abs{p}-1)!\prod_i(q_i+1)}\frac{(F_y)_{y,q}}{F_y^{\abs{p}+\abs{q}}}.
\end{equation} 
Substituting this into equation (\ref{LR1}) and simplifying gives
\begin{equation}
\lim_{y\rightarrow 0}(\frac{\partial^{\abs{p}-1}}{\partial y^{\abs{p}-1}}(\frac{F_{x,p}}{(F/y)^{\abs{p}}}))=\sum_{\substack{s=0\\q\vdash
s}}^{\abs{p}-1}\frac{\alpha_{q}(-1)^{\abs{q}}(\abs{p}+\abs{q}-1)!}{s!\prod{(q_i+1)}(\abs{p}-1-s)!}\frac{(F_y)_{y,q}}{F_y^{\abs{p}+\abs{q}}}\frac{\partial^{\abs{p}-1-s}F_{x,p}}{\partial y^{\abs{p}-1-s}}.
\end{equation} 
Then by equations (\ref{GLI2}) and (\ref{LR0}),
\begin{equation}\label{EF2}
\frac{d^ny}{dx^n}=\sum_{\substack {p\vdash n\\q\vdash s\leq\abs{p}-1}}\frac{(-1)^{\abs{p}+\abs{q}}
{\alpha_p\alpha_q(\abs{p}+\abs{q}-1)!}}{s!\prod{(q_i+1)}(\abs{p}-1-s)!}\frac{(F_y)_{y,q}}{F_y^{\abs{p}+\abs{q}}}
\frac{\partial^{\abs{p}-1-s}F_{x,p}}{\partial y^{\abs{p}-1-s}}.
\end{equation}
Although we were assuming $x=y=0$ in the derivation of equation (\ref{EF2}), the same equation clearly holds at any point. We use Leibniz' rule to write out
\begin{equation}
\frac{\partial^{\abs{p}-1-s}F_{x,p}}{\partial y^{\abs{p}-1-s}}=\frac{\partial^{\abs{p}-1-s}}{\partial y^{\abs{p}-1-s}}
(\prod \frac{\partial^{p_i}F}{\partial
x^{p_i}})=\sum_{\beta_1+...+\beta_{\abs{p}}=\abs{p}-1-s}\frac{(\abs{p}-1-s)!}{\beta_1!...\beta_{\abs{p}}!}\prod_{i=1}^{\abs{p}}\frac{\partial^{p_i+\beta_i}F}{\partial
x^{p_i}\partial y^{\beta_i}}.
\end{equation} 
Then equation (\ref{EF2}) becomes
\begin{equation}\label{EF3}
\frac{d^ny}{dx^n}=\sum_{\substack{p\vdash n\\q\vdash s\leq\abs{p}-1\\\beta_1,...,\beta_{\abs{p}}:\sum\beta_i=\abs{p}-1-s}}
\frac{(-1)^{\abs{p}+\abs{q}}\alpha_p\alpha_q(\abs{p}+\abs{q}-1)!}{\beta_1!...\beta_{\abs{p}}!s!\prod{(q_i+1)}}\frac{F_{x^{p_1}y^{\beta_1}}...F_{x^{p_{\abs{p}}}y^{\beta_{\abs{p}}}}(F_y)_{y,q}}{F_y^{\abs{p}+\abs{q}}}.
\end{equation}
In equation (\ref{EF3}), each term corresponds to a partition $p\vdash n,$ a partition $q\vdash s$ for some $s\leq\abs p-1,$ and $\abs{p}$ non-negative integers $\beta_1,...,\beta_{\abs p}$ with $\sum\beta_i=\abs p-1-s.$ For a given $p,$ $q$ and $(\beta_i),$ there is a unique permutation $\pi$ of $1,..,\abs p$ such that the resulting sequence $(p_{\pi(1)},\beta_{\pi(1)}),...,(p_{\pi(\abs p)},\beta_{\pi(\abs p)})$ is lexicographically ordered, and $p,\beta$ and $q$ thus can be combined in a unique way into a two-dimensional partition $P(p,\beta,q)$ given by
\begin{equation}\label{formP}
P(p,\beta,q):\overbrace{(p_{\pi(1)},\beta_{\pi(1)}), ... ,(p_{\pi(\abs p)},\beta_{\pi(\abs p)})}^{\text{lexicographically ordered}},(0,q_1+1),...,(0,q_{\abs q}+1).
\end{equation}
We see that $P=P(p,\beta,q)$ is lexicographically ordered throughout; indeed $p_{\pi(\abs p)}>0,$ so $(p_{\pi(\abs p)},\beta_{\pi(\abs p)})>(0,q_1+1).$ Also $P$ satisfies $\abs P=\abs p +\abs q,$ $P\vdash (n,\abs P-1)$ and $(0,1)\notin P.$ Furthermore, clearly any partition satisfying these two relations, is of the form (\ref{formP}). Moreover,
\begin{equation}
\frac{F_{x^{p_1}y^{\beta_1}}...F_{x^{p_{\abs p}}y^{\beta_{\abs{p}}}}(F_y)_{y,q}}{F_y^{\abs p+\abs q}}=\frac{F_{x,y,P}}{F_y^{\abs P}}.
\end{equation} 
In particular, this term is the same for each triple $(p,(\beta_i),q)$ in equation (\ref{EF3}) which gives rise to $P,$ so we can change the summation to be over partitions $P,$ which will lead us to the final form required. To do this, we need to calculate how many triples $(p,\beta_i,q)$ give rise to the same $P.$ But $P$ occurs as $P=P(p,(\beta_i),q)$ for a unique $p$ and $q,$ and for permutations of the $(\beta_i)$ which preserve the set \begin{equation}
\set{(p_1,\beta_1),...,(p_{\abs p},\beta_{\abs p})},
\end{equation}
so that after lexicographic ordering, the same partition $P$ is obtained. The number of such permutations is \begin{equation}\prod_k e_{p,k}!/\prod_{k>0,l\geq 0}e_{P,k,l}!\end{equation} (The product in the denominator is over $k>0,$ as the multiplicities $e_{P,0,l}$ are not relevant, being those associated with the part of $P$ coming from $q.)$ All the associated terms in equation (\ref{EF3}) are the same, so equation (\ref{EF3}) becomes
\begin{equation}\label{EF4}
\frac{d^ny}{dx^n}=\sum_{\substack{P\vdash (n,\abs{P}-1)\\(0,1)\notin P}}\frac{(-1)^{\abs{p}+\abs{q}}\alpha_p\alpha_q(\abs{p}+\abs{q}-1)!}{\beta_1!...\beta_{\abs{p}}!s!\prod_i{(q_i+1)}}\frac{\prod_k e_{p,k}!}{\prod_{k>0,l\geq 0}e_{P,k,l}!}\frac{F_{x,y,P}}{F_y^{\abs p}}.
\end{equation}
The coefficient in equation (\ref{EF4}) is 
\begin{equation}\label{EF3a}
\begin{split}
\frac{(-1)^{\abs{p}+\abs{q}}\alpha_p\alpha_q(\abs{p}+\abs{q}-1)!\prod_k e_{p,k}!}{\beta_1!...\beta_{\abs{p}}!s!\prod_i{(q_i+1)}\prod_{k>0,l\geq 0}e_{P,k,l}!}
&=\frac{(-1)^{\abs{p}+\abs{q}}n!s!(\abs{p}+\abs{q}-1)!}{\prod_i p_i!\prod_k e_{p,k}!\prod_i\beta_i!\prod_j
e_{qj}!\prod_i{(q_i+1)!}s!}\frac{\prod_k e_{p,k}!}{\prod_{k>0,l\geq 0}e_{P,k,l}!}\\
&=\frac{(-1)^{\abs{P}}n!(\abs{P}-1)!}{\prod_iP_{i1}!\prod_i{P_{i2}!}\prod_j
e_{P,0,j}!\prod_{k>0,l\geq 0}e_{P,k,l}!}\\
&=(-1)^{\abs P}\alpha_P,
\end{split}
\end{equation}
where we have used that $\prod P_{i1}!=\prod p_i!$ and $\prod P_{i2}!=\prod\beta_i!\prod (q_i+1)!,$ which are evident from equation (\ref{formP}) and the definitions of $\alpha_p$ from equations (\ref{alpha1}) and (\ref{alpha2}). Thus equation (\ref{EF4}) amounts to equation (\ref{main}), and the derivation is complete. 

\section{Proof of equation (\ref{main}) by induction}\label{S5}
It is possible to prove equation (\ref{main}) by induction on $n.$ The inductive proof is completely elementary and it therefore seems worthwhile to give this proof here. To simplify the presentation, we will find it useful to introduce the following operations on partitions. For a partition $p\vdash (n,m),$ and for $i,j$ with $i>0,$ $j\geq 0,$ $(i,j)\neq(1,0)$ and $e_{p,i,j}>0,$  define $q=A_{i,j}(p)$ to be the partition with
\begin{equation}\label{Aij}
e_{q,k,l}=\begin{cases}e_{p,k,l}-1 & (k,l)=(i,j)\\e_{p,k,l}+1 & (k,l)=(i-1,j)\\e_{p,k,l} &\text{ otherwise}\end{cases}.
\end{equation}
In words, viewing the table of muliplicities for $p,$ we move one unit from the $(i,j)$ position to the $(i-1,j)$ position. Notice the condition $(i,j)\neq(1,0)$ ensures that the operation results in a partition, i.e. we have $e_{A_{i,j}(p),0,0}=0.$ Similarly, for $i,j$ with $i\geq 0,$ $j>0,$ $(i,j)\neq(0,1)$ and $e_{p,i,j}>0,$  define $q=B_{i,j}(p)$ to be the partition with
\begin{equation}\label{Bij}
e_{q,k,l}=\begin{cases}e_{p,k,l}-1 & (k,l)=(i,j)\\e_{p,k,l}+1 & (k,l)=(i,j-1)\\e_{p,k,l} &\text{ otherwise}\end{cases}.
\end{equation}
The operations $A_{i,j}$ and $B_{i,j}$ are convenient for expressing the derivatives of a term $F_{x,y,p},$ via the following two formulae.
\begin{equation}\label{ddx}
\begin{aligned}
\frac{\partial}{\partial x}(F_{x,y,q})
&=\frac{\partial}{\partial x}(\prod_{i,j:e_{q,i,j}>0}(\frac{\partial^{i+j}F}{\partial x^i\partial y^j})^{e_{q,i,j}})\\
&=\sum_{i,j:e_{q,i,j}>0}\prod_{(k,l):k\neq i,l\neq j}(\frac{\partial^{k+l}F}{\partial x^k\partial y^l})^{e_{q,k,l}}e_{q,i,j}(\frac{\partial^{i+j}F}{\partial x^i\partial y^j})^{e_{q,i,j}-1}\frac{\partial^{i+j+1}F}{\partial x^{i+1}\partial y^j}\\
&=\sum_{i,j:e_{q,i,j}>0}e_{q,i,j}F_{x,y,A_{i+1,j}^{-1}(q)}\\
&=\sum_{i,j:e_{q,i-1,j}>0}e_{q,i-1,j}F_{x,y,A_{i,j}^{-1}(q)}\\
&=\sum_{p,i,j:A_{i,j}(p)=q}(e_{p,i-1,j}+1)F_{x,y,p},
\end{aligned}\end{equation}
where the sum is over all $(p,i,j)$ such that $A_{i,j}(p)$ is defined and equal to $q.$ Similarly,
\begin{equation}\label{ddy}
\frac{\partial}{\partial y}(F_{x,y,q})=\sum_{p,i,j:B_{i,j}(p)=q}(e_{p,i,j-1}+1)F_{x,y,p}.
\end{equation}
Suppose $p\vdash (n,m).$ Then it is easy to verify that
\begin{equation}\label{abspabsq}
\abs{A_{i,j}(p)}=\abs{B_{i,j}(p)}=\abs{p}
\end{equation}
and 
\begin{equation}\label{sumpsumq}
A_{i,j}(p)\vdash (n-1,m)\text{ and }B_{i,j}(p)\vdash (n,m-1),
\end{equation}
and further, using the definition at equation (\ref{alpha2}),
\begin{equation}\label{apaq1}
\alpha_{A_{i,j}(p)}=\frac{ie_{p,i,j}}{n(e_{p,i-1,j}+1)}\alpha_p\text{ and }\alpha_{B_{i,j}(p)}=\frac{je_{p,i,j}}{m(e_{p,i,j-1}+1)}\alpha_p.
\end{equation}

We also need the following, where the backslash notation indicates removing the specified parts from a partition $p\vdash(n,m)$ containing at least one of each part:
\begin{equation}\label{apaq2}
\alpha_{p\backslash (1,0)}=\frac{e_{p,1,0}}{n}\alpha_p,\alpha_{p\backslash (1,1)}=\frac{e_{p,1,1}}{nm}\alpha_p
\text{ and }
\alpha_{p\backslash (1,0),(0,2)}=\frac{2e_{p,1,0}e_{p,0,2}}{nm(m-1)}\alpha_p.
\end{equation}
Combining the second equation of (\ref{apaq1}) and the first equation of (\ref{apaq2}) gives, where again $p\vdash(n,m),$
\begin{equation}\label{apaq3}
\alpha_{B_{i,j}(p\backslash(1,0))}=\frac{je_{p\backslash(1,0),i,j}}{m(e_{p\backslash(1,0),i,j-1}+1)}\alpha_{p\backslash(1,0)}
=\frac{je_{p,i,j}e_{p,1,0}}{mn(e_{p\backslash(1,0),i,j-1}+1)}\alpha_p,
\end{equation}
wherever $B_{i,j}(p\backslash(1,0))$ is defined. (In the right hand we have replaced $e_{p\backslash(1,0),i,j}$ with $e_{p,i,j},$ since they are equal for $j>0,$ and $B_{i,j}$ is only defined for $j>0.$)

We are now ready to begin the proof of equation (\ref{main}) by induction. For $n=1,$ equation (\ref{main}) is just equation (1), so by induction assume equation (\ref{main}) holds for $n-1.$ Then by equation (\ref{TD}), we have 
\begin{equation}\label{applyTD}
\frac{d^ny}{dx^n}=\sum_{\substack{q\vdash (n-1,\abs{q}-1)\\(0,1)\notin q}}(-1)^{\abs q}\alpha_q(\frac{\partial}{\partial x}-\frac{F_x}{F_y}\frac{\partial}{\partial y})(F_y^{-\abs q}F_{x,y,q}),
\end{equation}
For any partition $q$ in the sum, we have four terms
\begin{equation}\label{X14}
(\frac{\partial}{\partial x}-\frac{F_x}{F_y}\frac{\partial}{\partial y})(F_y^{-\abs q}F_{x,y,q})=X_1+X_2+X_3+X_4,
\end{equation}
where
\begin{equation}\label{Xieqn}\begin{aligned}
X_1&=\frac{\partial}{\partial x}(F_y^{-\abs q})F_{x,y,q}=-\abs q F_y^{-\abs q-1}F_{x,y,q\cup(1,1)};\\
X_2&=F_y^{-\abs q}\frac{\partial}{\partial x}(F_{x,y,q})=F_y^{-\abs q}\sum_{p,i,j:A_{i,j}(p)=q}(e_{p,i-1,j}+1)F_{x,y,p};\\
X_3&=-\frac{F_x}{F_y}\frac{\partial}{\partial y}(F_y^{-\abs q})F_{x,y,q}=\abs q F_y^{-\abs q-2}F_{x,y,q\cup(1,0)\cup(0,2)};\\
X_4&=-F_y^{-\abs q-1}F_x\frac{\partial}{\partial y}(F_{x,y,q})=-F_y^{-\abs q-1}\sum_{p,i,j,B_{i,j}(p\backslash(1,0))=q}(e_{p\backslash(1,0),i,j-1}+1)F_{x,y,p},
\end{aligned}\end{equation}
and we have used equations (\ref{ddx}) and (\ref{ddy}) respectively, to obtain the expressions for $X_2$ and $X_4.$ The sums in $X_2$ and $X_4$ are thus over all triples $(p,i,j)$ where $p$ is a partition and $i,j$ are integers, satisfying that $A_{i,j}(p)=q$ respectively $B_{i,j}(p\backslash (1,0))=q.$

In equation (\ref{Xieqn}), we see four types of relationship between partitions $p$ and $q,$ corresponding to the four terms $X_{1-4}.$ The following lemma collates information about partitions that are related in one of these ways.

\begin{lemX}
Let $p$ and $q$ be two-dimensional partitions and let $p\vdash (n,m).$ Suppose that $p$ and $q$ are related by one of the following four equations.
\renewcommand\labelenumi{\theenumi.}
\begin{enumerate}
\item $q=p\backslash(1,1).$
\item $q=A_{i,j}(p)$ for some $i,j.$ 
\item $q=p\backslash(1,0),(0,2).$
\item $q=B_{i,j}(p\backslash (1,0))$ for some $i,j.$
\end{enumerate}
\renewcommand\labelenumi{(\theenumi)}
Then the following hold:
\begin{enumerate}
\item We have 
$$
\abs p=\begin{cases}\abs q+1&\text{in case }1\\\abs q&\text{in case }2\\\abs q+2&\text{in case }3\\\abs q+1&\text{in case }4\end{cases}\text{ and }q\vdash\begin{cases}(n-1,m-1)&\text{in case }1\\(n-1,m)&\text{in case }2\\(n-1,m-2)&\text{in case }3\\(n-1,m-1)&\text{in case }4\end{cases}.$$
\item $p\vdash (n,\abs{p}-1)$ if and only if $q\vdash (n-1,\abs{q}-1).$
\item If $(0,1)\notin q$ then $(0,1)\notin p,$ and if $(0,1)\notin p$ then $(0,1)\notin q$ unless either $p$ and $q$ are related by case $2$ with $(i,j)=(1,1),$ or $p$ and $q$ are related by case $4$ with $(i,j)=(0,2).$ In either of these exceptions, we always have $(0,1)\in q.$
\item $$\alpha_q=\begin{cases}\frac{e_{p,1,1}}{nm}\alpha_p&\text{in case }1\\\frac{ie_{p,i,j}}{n(e_{p,i-1,j}+1)}\alpha_p&\text{in case }2\\\frac{2e_{p,1,0}e_{p,0,2}}{nm(m-1)}\alpha_p&\text{in case }3\\\frac{je_{p,i,j}e_{p,1,0}}{nm(e_{p\backslash(1,0),i,j-1}+1)}\alpha_p&\text{in case }4\end{cases}.$$
\end{enumerate}
\end{lemX}
\begin{proof}
(1) The assertions about $\abs p$ are clear in cases $1$ and $3,$ while cases $2$ and $4$ are immediate from equations (\ref{abspabsq}). Similarly, in cases $1$ and $3$ it is clear that $q\vdash (n-1,m-1)$ and $q\vdash (n-1,m-2)$ respectively, while equations (\ref{sumpsumq}) show that $p\vdash (n-1,m)$ and $p\vdash (n-1,m-1)$ in cases $2$ and $4$ respectively. This gives (1). 

(2) follows directly from (1).

(3) First suppose $(0,1)\notin q.$ Then it is clear in cases $1$ and $3$ that $(0,1)\notin p.$ For case $2,$ we need to show that $(0,1)\in p$ implies $(0,1)\in A_{i,j}(p)$ whenever the latter is defined. This follows from the definition of $A_{i,j}$ at equation (\ref{Aij}). Similarly, for case $4,$ we need that $(0,1)\in p$ implies $(0,1)\in B_{i,j}(p)$ whenever $B_{i,j}(p)$ is defined, which again follows from the definition (\ref{Bij}). 

Conversely, suppose $(0,1)\notin p.$ Then again it is obvious in cases $1$ and $3$ that $(0,1)\notin q.$ In cases $2$ and $4,$ we have from the definitions (\ref{Aij}) and (\ref{Bij}), that $(0,1)\notin p$ implies $(0,1)\notin A_{i,j}(p),$ unless $(i,j)=(1,1),$ when $(0,1)\in A_{i,j}(p),$ and similarly $(0,1)\notin p$ implies $(0,1)\notin B_{i,j}(p),$ unless $(i,j)=(0,2),$ when $(0,1)\in B_{i,j}(p).$ 

(4) merely collates the second equation in (\ref{apaq2}), the first equation in (\ref{apaq1}), the third equation in (\ref{apaq2}) and the last equation in (\ref{apaq3}), for cases $1-4$ respectively.
\end{proof}

We now complete the proof by induction. First, by comparing the powers of $F_y$ appearing in each term $X_{1-4}$ in equation (\ref{X14}), with the assertions about $\abs p$ in part (1) of the lemma, we see that each term is of the form $F_y^{-\abs p}F_{x,y,p},$ where $p$ and $q$ are related by one of the cases $1-4$ of the lemma. By parts (2) and (3) of the lemma, since $q\vdash (n-1,\abs q-1)$ and $(0,1)\notin q,$ we have $p\vdash (n,\abs p-1)$ and $(0,1)\notin p.$ We have therefore established that equation (\ref{applyTD}) is of the form
\begin{equation}\label{tpeqn}
\frac{d^ny}{dx^n}=\sum_{\substack{p\vdash (n,\abs{p}-1)\\(0,1)\notin p}}t_p\frac{F_{x,y,p}}{F_y^{\abs p}}.
\end{equation}
for some coefficients $t_p.$ It remains to show that $t_p=(-1)^{\abs p}\alpha_p.$ To do this, we fix a particular $p\vdash (n,m)$ where $m=\abs p-1$ and $(0,1)\notin p.$ By equations (\ref{applyTD}) and (\ref{X14}), we have
\begin{equation}\label{ABCD}
t_p=A+B+C+D,
\end{equation}
where $A,B,C,D$ are respectively the contributions to $t_p$ from all $q\vdash (n-1,\abs q-1)$ with $(0,1)\notin q,$ via the terms $X_{1-4}$ in equation (\ref{Xieqn}). These are given by:
\begin{equation}\label{A-D}
\begin{aligned}
A&=-(-1)^{\abs q}\alpha_q\abs q \text{ if }(1,1)\in p\text{ and }A=0\text{ otherwise,}\\
B&=\sum_{q,i,j:A_{i,j}(p)=q:(i,j)\neq(0,2)}(e_{p,i-1,j}+1)(-1)^{\abs q}\alpha_q,\\
C&=(-1)^{\abs q}\alpha_q\abs q\text{ if }(1,0),(0,2)\in p\text{ and }C=0\text{ otherwise,}\\
D&=-\sum_{q,i,j:B_{i,j}(p\backslash(1,0))=q:(i,j)\neq(1,1)}(e_{p\backslash(1,0),i,j-1}+1)(-1)^{\abs q}\alpha_q\abs q.
\end{aligned}
\end{equation}
Note that in equations (\ref{A-D}), we have the conditions $(1,1)\notin p$ and $(0,2)\notin p$ respectively in the sums for $B$ and $D.$ By parts (2) and (3) of the lemma, these are the conditions needed to ensure that $(0,1)\notin q$ and $q\vdash (n-1,\abs q-1).$

Using parts (2) and (4) of the lemma, we compute
\begin{equation}\label{A}
A=-(-1)^{\abs q}\abs q\frac{e_{p,1,1}}{nm}\alpha_p=(-1)^{\abs p}(\abs p-1)\frac{e_{p,1,1}}{nm}\alpha_p=(-1)^{\abs p}\frac{e_{p,1,1}}{n}\alpha_p,
\end{equation}
noting that this equation automatically gives the correct zero contribution when $(1,1)\notin p.$ Next,
\begin{equation}\label{B}
B=\sum_{\substack{q,i,j:A_{i,j}(p)=q\\(i,j)\neq(1,1)}}(e_{p,i-1,j}+1)(-1)^{\abs p}\frac{ie_{p,i,j}}{n(e_{p,i-1,j}+1)}\alpha_p=\sum_{(i,j)\neq(1,0),(1,1)}(-1)^{\abs p}\frac{ie_{p,i,j}}{n}\alpha_p.
\end{equation}
Note that the range of summation in equation (\ref{B}) is the set $(i,j)$ having $A_{i,j}(p)$ defined and $(i,j)\neq(1,1).$ By equation (\ref{Aij}), $A_{i,j}(p)$ is defined where $i>0,j\geq 0$ $e_{p,i,j}>0$ and $(i,j)\neq(1,0).$ However if $i=0$ or $e_{p,i,j}=0$ then the corresponding summand is automatically zero, so the conditions $i>0$ and $e_{p,i,j}>0$ are removed in the right hand expression.

Similarly to equation (\ref{A}), we have
\begin{equation}\label{C}
C=(-1)^{\abs q}\frac{2e_{1,0}e_{0,2}}{nm(m-1)}\alpha_p\abs q=(-1)^{\abs p}\frac{2e_{p,1,0}e_{p,0,2}}{nm}\alpha_p,
\end{equation}
and similarly to equation (\ref{B}),
\begin{equation}
D=\sum_{\substack{(i,j),i\geq 0,j>0,\\(i,j)\neq(1,0),(0,2)}}-(e_{p\backslash(1,0),i,j-1}+1)\frac{(-1)^{\abs q}e_{p,1,0}je_{p,i,j}}{nm(e_{p\backslash(1,0),i,j-1}+1)}\alpha_p\abs q=\sum_{(i,j)\neq(0,2)}(-1)^{\abs p}\frac{e_{p,1,0}je_{p,i,j}}{n}\alpha_p.
\end{equation}
Substituting in equation (\ref{ABCD}) gives
\begin{equation}
\begin{aligned}
t_p&=(-1)^{\abs p}\alpha_p\left\{\frac{e_{p,1,1}}{n}+\sum_{(i,j)\neq(1,0),(1,1)}\frac{ie_{p,i,j}}{n}+\frac{2e_{p,1,0}e_{p,0,2}}{nm}+\sum_{(i,j)\neq(0,2)}\frac{je_{p,1,0}e_{p,i,j}}{nm}\right\}\\
&=(-1)^{\abs p}\alpha_p\left\{\sum_{(i,j)\neq(1,0)}\frac{ie_{p,i,j}}{n}+\sum_{i,j}\frac{je_{p,1,0}e_{p,i,j}}{nm}\right\}\\
&=(-1)^{\abs p}\alpha_p\left\{1-\frac{e_{p,1,0}}{n}+\frac{e_{p,1,0}}{n}\right\}\\
&=(-1)^{\abs p}\alpha_p,
\end{aligned}
\end{equation}
where we have used $\sum_{i,j}ie_{p,i,j}=n$ and $\sum_{i,j}je_{p,i,j}=m.$ This completes the proof.

\section{Discussion of Comtet and Fiolet's formula}\label{S6}
In this section, we compare Comtet and Fiolet's formula in \cite{Comtet1974} with our equation (\ref{main}). Comtet and Fiolet write $E=\mathbb N\times\mathbb N\backslash\set{(0,0),(0,1)}$ where $\mathbb N$ denotes the natural numbers including zero. Also, in a notation slightly different from ours above, they let $F_{i,j}$ denote $\partial^{i+j}F/\partial x^i\partial y^j.$ Then their formula (\cite[Th\'eor\`eme 1]{Comtet1974}) is 
\begin{equation}\label{CF1}
\frac{d^ny}{dx^n}=\sum_{m=1}^{2n-1}(\frac{-1}{F_y})^mI_{n,m}
\end{equation}
where
\begin{equation}\label{CF2}
I_{n,m}=\sum_{\substack{l_1+2l_2+3l_3+...=n\\c_1+2c_2+3c_3+...=m-1}}\frac{n!q!\langle q\rangle_S\langle q+S\rangle_{c_1}}{\prod_{k\geq 1}(k!)^{c_k+l_k}}\prod_{(i,j)\in E}\frac{F_{i,j}^{d_{i,j}}}{d_{i,j}!},
\end{equation}
the sum being over tables of non-negative integers $(d_{i,j}):(i,j)\in E,$ and where the following notation is used:
\begin{equation}\label{CF3}
\begin{aligned}
l_i=\sum_jd_{i,j};\ c_j&=\sum_id_{i,j}\text{ (row and column sums)},\\
S=\sum_{j\geq 2}c_j;\ q&=1+\sum_{j\geq 1}jc_{j+1},\\
\langle q\rangle_S&=q(q+1)...(q+S-1).
\end{aligned}
\end{equation}

To see the basic correspondence between equations (\ref{CF2}) and (\ref{main}), we regard each table of non-negative integers $(d_{i,j}):(i,j)\in E$ in equation (\ref{CF2}), as the multiplicities of a unique two-dimensional partition $p,$ so that $d_{i,j}=e_{p,i,j}.$ We have $(0,1)\notin p,$ by definition of $E.$ The conditions under the sum in equation (\ref{CF2}) are then equivalent to requiring that $p\vdash (n,m-1).$ To assist comparison, we rewrite our equation (\ref{main}), to look similar to equations (\ref{CF1}) and (\ref{CF2}), putting $m=\abs p$ to give
\begin{equation}\label{main1}
\frac{d^ny}{dx^n}=\sum_{m=1}^{2n-1}(\frac{-1}{F_y})^m\sum_{\substack{p\vdash (n,m-1)\\ \abs p=m\\(0,1)\notin p}}\alpha_pF_{x,y,p}.
\end{equation} 
The upper limit of summation $m=2n-1$ in equation (\ref{main1}) is automatic, since if $p\vdash(n,\abs p-1)$ and $(0,1)\notin p$ as required in our formula, then $\abs p\leq 2n-1.$ To see this, write $p$ in the form \begin{equation}(p_1,q_1),...,(p_r,q_r),(0,q_{r+1}),...,(0,q_{r+s}),\end{equation} with $p_r>0.$ Then $\abs p=r+s,$ and since $q_{r+1},...,q_{r+s}\geq 2,$ we have $r+s-1=\abs p-1\geq 2s,$ so $r-1\geq s,$ while also $n\geq r,$ so $2n-1\geq 2r-1\geq r+s=\abs p.$ 

Now for a partition $p,$ to assume Comtet and Fiolet's notation, write $d_{i,j}$ instead of $e_{p,i,j}.$ Then
\begin{equation}
\prod_ip_{i1}!=\prod_k (k!)^{c_k}\text{ and }\prod_ip_{i2}!=\prod_k (k!)^{l_k},\end{equation}
where $c_k$ and $l_k$ are defined in equation (\ref{CF3}). Thus since $p\vdash (n,m-1),$ we may write equation (\ref{main1}) as
\begin{equation}\label{likeCF}
\frac{d^ny}{dx^n}=\sum_{m=1}^{2n-1}(\frac{-1}{F_y})^m\sum_{\substack{l_1+2l_2+3l_3+...=n\\c_1+2c_2+3c_3+...=m-1\\l_0+l_1+l_2...=m}}\frac{n!(m-1)!}{\prod_k (k!)^{c_k+l_k}}\prod_{i,j\in E}\frac{F_{i,j}^{d_{i,j}}}{d_{i,j}!}.
\end{equation}
There are just two differences between our equation (\ref{likeCF}) and Comtet and Fiolet's formula comprised in (\ref{CF1}) and (\ref{CF2}).
\begin{enumerate}
\item We have the additional condition $l_0+l_1+..=\sum_{i,j}d_{i,j}=m$ under the inner sum. 
\item We have a term $(m-1)!$ in the coefficient, where Comtet and Fiolet have the term $q!\langle q\rangle_S\langle q+S\rangle_{c_1}.$
\end{enumerate} 
We deal with these points in turn. $(1)$ Our additional condition $l_0+l_1+..=\sum_{i,j}d_{i,j}=m$ under the inner sum, ensures that $\abs p=m,$ or in other words that $p\vdash (n,\abs p-1).$ There is no equivalent condition in equation (\ref{CF2}), but in the text immediately following their statement, Comtet and Fiolet observe that the polynomial $I_{n,m}$ defined by equation (\ref{CF2}) is homogeneous of degree $m,$ in the sense that $\sum_{i,j}d_{i,j}=m.$ This is not an automatic consequence of the other conditions: Indeed, equation (\ref{CF2}) as it stands includes terms for any partition $p\vdash (n,m-1)$ with $m\leq 2n-1,$ and these certainly do not all satisfy $\abs{p}=m.$ (Perhaps the extra condition was simply lost during the typesetting of (\ref{CF2}); there is a row of dots under the condition $c_1+2c_2+3c_3+...=m-1,$ in equation (\ref{CF2}) as it appears in \cite{Comtet1974}, which look like separators for a further line below that was not printed.)

Turning to difference $(2),$ Comtet and Fiolet have a term $q!\langle q\rangle_S\langle q+S\rangle_{c_1},$ where in our equation (\ref{likeCF}) we just have $(m-1)!.$ However, these terms are nearly equal (but our version is correct, as we check below with an example). From the definitions at equation (\ref{CF3}) we have
\begin{equation}
q+S+c_1=1+\sum_{j\geq 1}jc_{j+1}+\sum_{j\geq 2}c_j+c_1=1+\sum_{j\geq 1}jc_j=m,
\end{equation}
and hence in equation (\ref{CF2}),
\begin{equation}
q!\langle q\rangle_S\langle q+S\rangle_{c_1}=q!(q+1)...(q+S)(q+S+1)...(q+S+c_1-1)=q(m-1)!.
\end{equation}
Thus, equation (\ref{CF2}) becomes 
\begin{equation}\label{CF4}
I_{n,m}=\sum_{\substack{l_1+2l_2+3l_3+...=n\\c_1+2c_2+3c_3+...=m-1\\l_0+l_1+l_2...=m}}\frac{n!q(m-1)!}{\prod_{k\geq 1}(k!)^{c_k+l_k}}\prod_{(i,j)\in E}\frac{F_{i,j}^{d_{i,j}}}{d_{i,j}!},
\end{equation}
Now, the terms of equations (\ref{CF4}) and (\ref{likeCF}) agree except that there is an additional factor of $q$ in (\ref{CF4}). To check that our version (\ref{likeCF}) is indeed the correct version, consider the case $n=5,$ and the coefficient of the term belonging to the partition $p\vdash (5,4)$ given by
\begin{equation}
(5,4)=(1,1)+(1,1)+(1,1)+(1,0)+(1,0)+(0,2).
\end{equation}
Note that this partition satisfies $p\vdash (5,\abs p-1)$ and $(0,1)\notin p.$ The coefficient of the corresponding term $F_x^2F_{xy}^3F_{yy}$ in the expression for $d^5y/dx^5$ is $600.$ This is stated correctly by Comtet and Fiolet in the numerical table in \cite{Comtet1974} (it is the coefficient of $f_{1,0}^2f_{1,1}^3f_{0,2}$ in their notation), and we have confirmed its value by performing a recursive calculation up to $n=6$ using equation (\ref{TD}). However for this partition, we have 
\begin{equation}
d_{1,0}=2,d_{1,1}=3,d_{0,2}=1,
\end{equation}
and all other $d_{i,j}$ are zero. Hence the terms going into Comtet and Fiolet's formula are as follows:
\begin{equation}
\begin{aligned}
n&=5,\\
l_0&=1,l_1=5,l_2=...=0,\\
c_0&=2,c_1=3,c_2=1,c_3=...=0,\\
S&=c_2+c_3+...=1,\\
q&=1+c_2+2c_3+...=2.
\end{aligned}
\end{equation}
aso the coefficient according to Comtet and Fiolet's formula (\ref{CF2}) is
\begin{equation}
\frac{5!2!\langle2\rangle_1\langle3\rangle_3}{(2!)^11!2!3!}=\frac{120\times 2\times 2\times 3\times 4\times 5}{4\times 6}=1200,
\end{equation}
which is out by a factor of $q=2,$ as expected, while our formula (\ref{main}) gives the coefficient correctly as 
\begin{equation}
\alpha_p=\frac{5!4!}{2!2!3!}=600.
\end{equation}
In summary, Comtet and Fiolet's formula, once corrected, should look like the equation below, which differs only notationally from equation (\ref{main}).
\begin{equation}\label{CFcorrected}
\frac{d^ny}{dx^n}=\sum_{m=1}^{2n-1}(\frac{-1}{F_y})^m\sum_{\substack{l_1+2l_2+3l_3+...=n\\c_1+2c_2+3c_3+...=m-1\\c_0+c_1+c_2+...=m}}\frac{n!(q-1)!\langle q\rangle_S\langle q+S\rangle_{c_1}}{\prod_{k\geq 1}(k!)^{c_k+l_k}}\prod_{(i,j)\in E}\frac{F_{i,j}^{d_{i,j}}}{d_{i,j}!},
\end{equation}
with the proviso that the expression $(q-1)!\langle q\rangle_S\langle q+S\rangle_{c_1}$ is just a complicated way of writing $(m-1)!$

\section{The number of terms in equation (\ref{main})}\label{S7}
Comtet and Fiolet gave a formula (\cite[Th\'eor\`eme 2]{Comtet1974}) for the number of distinct terms appearing in their formula (\ref{CF1}-\ref{CF2}), as a function of $n.$ In terms of our equation (\ref{main}), this is the number of partitions $p\vdash (n,\abs p-1)$ with $(0,1)\notin p.$ Let this number be $a(n).$ Comtet and Fiolet gave the following table of values of $a(n)$ for $n\leq 23$ (we have added $a(24)$ from our own calculations).

$$\begin{tabular}{|p{.6cm}|p{.7cm}|p{.7cm}|p{.7cm}|p{.8cm}|p{.8cm}|p{.8cm}|p{.8cm}|p{1cm}|p{1cm}|p{1cm}|p{1cm}|p{1.1cm}|}
\hline
$n$&1&2&3&4&5&6&7&8&9&10&11&12\\
\hline
$a(n)$&\footnotesize 1&\footnotesize 3&\footnotesize 9&\footnotesize 24&\footnotesize 61&\footnotesize 145&\footnotesize 333&\footnotesize 732&\footnotesize 1565&\footnotesize 3247&\footnotesize 6583&\footnotesize 13047\\
\hline
$n$&13&14&15&16&17&18&19&20&21&22&23&24\\
\hline
$a(n)$&\footnotesize{25379}&\footnotesize{48477}&\footnotesize{91159}&\footnotesize{168883}&\footnotesize{308736}&\footnotesize{557335}&\footnotesize{994638}&\footnotesize{1755909}&\footnotesize{3068960}&\footnotesize{5313318}&\footnotesize{9118049}&\footnotesize{15516710}\\
\hline
\end{tabular}$$
\\

The sequence $1, 3, 9, 24, 61, ...$ of values $a(n)$ can be found in the On-Line Encyclopedia of Integer Sequences \cite{OEIS}. The numbers $a(n)$ were correctly stated by Comtet and Fiolet, but their formula (\cite[Th\'eor\`eme 2]{Comtet1974}) is not correct. This formula states that $a(n)$ is the coefficient of $t^nu^{n-1}$ in $\prod_{(i,j)\in E}(1-t^iu^j)^{-1},$ where $E=\mathbb N\times\mathbb N\backslash\set{(0,0),(0,1)}$ as defined above equation (\ref{CF1}). This is seen to be incorrect, for example, by calculating the coefficient of $t^2u.$ Expanding $\prod_{(i,j)\in E}(1-t^iu^j)^{-1},$ to all the terms that can possibly contribute, we get 
\begin{equation}
(1+t+t^2)(1+tu)(1+t^2)(1+t^2u),
\end{equation} 
from which we see that the coefficient of $t^2u$ in $\prod_{(i,j)\in E}(1-t^iu^j)^{-1}$ is $2,$ while $a(2)=3$ from the table above (or from equation (\ref{EG})). Comtet and Fiolet's incorrect formula is repeated on page 175 of \cite{Comtet1974a}. The correct result is as follows.
\begin{thmX}[Compare with \cite{Comtet1974}, Th\'eor\`eme 2 or \cite{Comtet1974a}, p175]
The number $a(n)$ of partitions $p\vdash (n,\abs p-1)$ and $(0,1)\notin p,$ or equivalently, the number $a(n)$ of distinct terms in equation (\ref{main}), is given by:
\begin{equation}\label{genfun}
a(n)=\mathrm{Coefficient\  of\ }t^nu^{n-1}\ \mathrm{in}\ \prod_{(i,j)\in E}\frac{1}{1-t^iu^{i+j-1}}
\end{equation}
\end{thmX}
\begin{proof}
Write
\begin{equation}
F(t,u)=\prod_{(i,j)\in E}\frac{1}{1-t^iu^{i+j-1}}=\prod_{(i,j)\in E}(\sum_{d\geq 0}(t^iu^{i+j-1})^d)
\end{equation}
For a given partition $p\vdash (n,\abs p-1)$ with $(0,1)\notin p,$ let $d_{i,j}=e_{p,i,j};(i,j)\in E,$ be the table of multiplicities of $p.$  Corresponding to this $p,$ we pick the term with $d=d_{i,j}$ in the $(i,j)$ factor of $F,$ and get a term of $F$ equal to \begin{equation}
t^{\sum_{i,j}id_{i,j}}u^{\sum_{i,j}(i+j-1)d_{i,j}}.
\end{equation}
This term contributes $1$ to the coefficient of $t^nu^{n-1}$ exactly when $$\sum_iid_{i,j}=n\text{ and }\sum_{i,j}(i+j-1)d_{i,j}=n-1$$ or \begin{equation}
\sum_iid_{i,j}=n\text{ and }\sum_{i,j}jd_{i,j}=\sum_{i,j}d_{i,j}-1,
\end{equation}
that is, exactly when $p\vdash (n,\abs p-1)$ as required.
\end{proof}

The numbers $a(n)$ can be calculated easily from the generating function (\ref{genfun}), and we indicate the simple scheme we used to recalculate the figures shown in the table (which, as we mentioned above, agree with the values given by Comtet and Fiolet for $n\leq 23$). Write 
\begin{equation}
F(u,t)=\prod_{(i,j)\in E}\frac{1}{1-t^iu^{i+j-1}}=\sum_{n\geq 0}p_n(t)u^n,
\end{equation}
so that $a(n)$ is the coefficient of $t^n$ in $p_{n-1}(t).$ Also let 
\begin{equation}
G(u,t)=\mathrm{log} F(u,t)=\sum_{(i,j)\in E}-\mathrm{log}(1-t^iu^{i+j-1})=\sum_{(i,j)\in E,r>0}\frac{(t^iu^{i+j-1})^r}{r}=\sum_{m\geq 0}q_m(t)u^m.
\end{equation}
For $m>0,$ we see that 
\begin{equation}\label{qmt}
q_m(t)=\sum_{\substack{i,j\geq 0\\i+j-1\text{ divides }m}}t^{im/(i+j-1)}=\sum_{d \text{ divides }m}\sum_{i=0}^{m/d+1}\frac{t^{id}}{d},
\end{equation}
from which the $q_m(t)$ are easy to compute explicitly. Now $F=e^G,$ so differentiating we have $FG_u=F_u,$ or 
\begin{equation}
(\sum_{r\geq 0}p_r(t)u^r)(\sum_{s>0}sq_s(t)u^{s-1})=\sum_{n>0}p_n(t)nu^{n-1},
\end{equation}
and comparing the coefficient of $u^{n-1}$ gives
\begin{equation}\label{Recurrence}
p_n(t)=\frac{1}{n}\sum_{s=1}^nsq_sp_{n-s}.
\end{equation} 
Starting from $p_0(t)=F(0,t)=1/(1-t)=\sum_{n\geq 0}t^n,$ equations (\ref{qmt}) and (\ref{Recurrence}) enable successive calculation of the polynomials $p_n(t)$ for $n>0,$ and then $a(n)$ is extracted as the coefficient of $t^n$ in $p_{n-1}(t).$ If we desire the values of $a(n)$ for $1\leq n\leq E,$ say, then we only have to keep terms up to $t^E$ in the polynomials $p$ and $q.$ To allow the reader to verify that this calculation does produce the numbers $a(n)$ in the table above (and hence also the numbers stated by Comtet and Fiolet), below is some self-explanatory VBA code implementing the calculation of $a(n)$ for $n\leq 50,$ say.
\\

\scriptsize
\begin{verbatim}
Sub Calc_AofN_upto_E()
	E = 50
	
	ReDim p(0 To E - 1, 0 To E): ReDim q(0 To E - 1, 0 To E)
	
	For m = 1 To E - 1: For d = 1 To m
      If m = d * Int(m / d) Then
           For i = 0 To m / d + 1
                If d * i <= E Then q(m, i * d) = q(m, i * d) + 1 / d
	Next: End If: Next: Next
	
	For j = 0 To E
  		p(0, j) = 1
	Next
	
	For n = 1 To E - 1: For s = 0 To n: For j = 0 To E: For i = 0 To j
      p(n, j) = p(n, j) + 1 / n * s * q(s, j - i) * p(n - s, i)
	Next: Next: Next: Next
	
	For n = 1 To E
    	Debug.Print p(n - 1, n)
	Next
End Sub
\end{verbatim}

London, UK. Email: tom@beech84.fsnet.co.uk

\normalsize

\begin{thebibliography}{2}
\bibitem {AbSteg} Abramowitz and Stegun, \emph{A handbook of mathematical functions}.
\bibitem{Comtet1968} Comtet. L., \emph{Polyn\^omes de Bell et formule explicite des d\'eriv\'ees successives d'une fonction implicite}, C.R. Acad. Sc. Paris S\'erie A, Vol 267, p 457-460.
\bibitem{Comtet1974} Comtet. L and Fiolet, M., \emph{Sur les d\'eriv\'ees successives d'une fonction implicite}, C.R. Acad. Sc. Paris S\'erie A, Vol
278, p 249-251.
\bibitem{Comtet1974a} Comtet. L., \emph{Advanced Combinatorics}, D. Reidel, Dordrecht, 1974.
\bibitem{Johnson} Johnson, W.P., \emph{The curious history of Fa\`{a} di Bruno's formula}, American Mathematical Monthly, No 109, March 2002, p 217-234.
\bibitem{OEIS} The On-Line Encyclopedia of Integer Sequences, http:\\www.research.att.com/~njas/sequences/A003262.

\end{thebibliography}
\end{document}